\documentclass{article}
\usepackage[utf8]{inputenc}
\usepackage{amsmath,amsthm,amssymb,amsfonts}
\usepackage[usenames]{color}

\usepackage{graphicx}
\usepackage{todonotes}
\usepackage{tikz} 
\usepackage{comment}

\def\1{{\bf 1}}

\title{How far from the edge need a population be to survive? A probability model.
}
\begin{document}

\maketitle

Rinaldo B. Schinazi
\footnote{Department of Mathematics, University of Colorado, Colorado Springs, CO 80933-7150, USA.
E-mail: rinaldo.schinazi@uccs.edu}

\begin{abstract}
Let $N$ be a natural number. We consider a population which lives on $I_N=\{-N,-N+1,\dots,N-1,N\}$. Each individual gives birth at rate $\lambda$ on each of its neighboring sites and dies at rate 1.
No births are allowed from the inside of $I_N$ to the outside or vice-versa.
 The population on the whole line (i.e. $N=+\infty$) survives with positive probability if and only if $\lambda>1/2$. On the other hand for any $1/2< \lambda\leq \sqrt 2/2$ there exists a natural number $N_c$ such that the population survives on $I_N$ for $N\geq N_c$ but dies out for $N<N_c$. 
There is no limit on the number of individuals per site so the population could grow at the center where the birth rates are maximum. Our result shows that it does not if the edge is too close.

\end{abstract}

\section{The model} A branching random walk (brw in short) on the one-dimensional lattice $\mathbb{Z}$ evolves as follows.
\begin{itemize}
    \item Each individual dies at rate 1.
    \item Let $x$ and $y$ in $\mathbb{Z}$ such that $|x-y|=1$. An individual at $x$ gives birth to an individual at $y$ at rate $\lambda$.
    \item There is no limit on the number of individuals per site.
\end{itemize}

Let $N\geq 1$ be a natural number. We are interested in a branching random walk restricted to $I_N=\{-N,-N+1,\dots,N-1,N\}$ with the following boundary conditions. No births are allowed from the inside of $I_N$ to the outside or vice-versa. In other words, if the population is to survive it has to survive on its own on the finite set of sites $I_N$. 

\section{Results}

We will show that if $\lambda>\sqrt 2/2$ the population on $I_N$ survives for all $N\geq 1$. On the other hand if $\lambda\leq 1/2$ the population dies out for all $N$. 

When $1/2<\lambda\leq \sqrt 2/2$ things get more interesting. Our main result is the following.

\medskip

{\bf Theorem 1.} {\sl Let $1/2<\lambda\leq \sqrt 2/2$ then there exists a natural number $N_c$ such that if $N< N_c$ the population restricted to $I_N$ dies out with probability 1 while if $N\geq N_c$ the population has a positive probability of surviving (i.e there is at least one individual alive in $I_N$ at all times).}

\section{Discussion}

Our result may have some relevance from a theoretical ecology point of view. It is well known that habitat fragmentation is one of the main causes of species extinction, see for instance \cite{Haddad}. Fragmentation results in less contiguous space but also in lower quality of the habitat. As a result one expects lower birth rates and higher death rates. This in turn can cause extinction of a population. In contrast to this scenario, our model shows that even with no change in birth and death rates the population will die out for being too close to the edge (i.e. $N$ too small). Note that there is no limit on the number of individuals per site so the population could grow at the center where the birth rates are maximum. Our result shows that it does not if the edge is too close.
Moreover, in this model the boundary is minimum (two sites at $-N$ and $N$). This should be helpful to the population since boundary sites are the only ones with lower birth rate. But even with minimal boundary our model suggests that fragmentation can be fatal.

\section{Proof of Theorem 1}

\subsection{The construction}

We start by giving an informal construction of the process. We first construct the process on the whole line $\mathbb{Z}$. The same construction will be used to construct the process on $I_N$ for all $N\geq 1$. At time 0 the initial configuration is assumed to have finitely many individuals on finitely many sites of $\mathbb{Z}$.  Every individual is assigned two Poisson processes, each with rate $\lambda$.
If the individual is at $x\in\mathbb{Z}$ then at each occurrence of the first Poisson process a new particle is born at $x-1$. Similarly, at each occurrence of the second Poisson process a new particle is born at $x+1$. Every individual is also assigned an exponential random variable with rate 1. At this exponential (random) time, the individual dies. To each new individual we again assign an exponential rate 1 random variable and two rate $\lambda$ Poisson processes and so on. All exponential random variables and Poisson processes are independent.

Let $N\geq 1$, we use the same exponential random variables and Poisson processes defined above to construct the process on $I_N$. The only difference with the construction on the whole line is that we suppress births from inside to outside and from outside to inside the box $I_N$. This allows the construction of the process on $I_N$ for all $N$ on the same probability space. 

\subsection{Monotone properties}

We now use this construction to show that the process (restricted or unrestricted) is increasing in $\lambda$. Assume that $\lambda_1<\lambda_2$. We construct the process with birth rate $\lambda_2$ and death rate 1 as indicated above. From this construction we get the process with rate $\lambda_1$ by filtering the births. That is, every time the $\lambda_2$ process has a birth we flip a coin. With probability $\lambda_1/\lambda_2$ the birth happens for the $\lambda_1$ process. With probability $1-\lambda_1/\lambda_2$ the birth does not happen for the $\lambda_1$ process. This is a well-known procedure to obtain a rate $\lambda_1$ Poisson process from a rate $\lambda_2$ Poisson process (see Section 3 in Chapter 12 in \cite{Schinazi}, for instance). This coupling shows that at every fixed time and for every site of $\mathbb{Z}$ the $\lambda_2$ process has more individuals than the $\lambda_1$ process. This is what we mean by stating that the process is increasing in $\lambda$.

We now turn to a monotone property in $N$. Let $N_1<N_2$, two natural numbers. By the construction above we can simultaneously construct the processes on $\mathbb{Z}$, restricted to $I_{N_1}$ and restricted $I_{N_2}$. To get the process restricted to $I_{N_2}$ from the one on $\mathbb{Z}$ we suppress the births from inside $I_{N_2}$ to outside $I_{N_2}$ and vice versa. We do the same for the process restricted to $I_{N_1}$. Since $N_1<N_2$, if a birth does not happen for the process restricted to $I_{N_2}$ it certainly does not happen for the process restricted to $I_{N_1}$. However, there are births that happen in $I_{N_2}$ but not in $I_{N_1}$. Hence, at every fixed time and for every site the process restricted to $I_{N_2}$ has more individuals  than the process restricted to $I_{N_1}$. In this sense the process is increasing in $N$.

\subsection{Tom Liggett's result}

We first deal with the case $\lambda<1/2$. 
Consider the process on the whole line $\mathbb{Z}$. Every individual has a birth rate of $2\lambda$ and a death rate of 1. Hence, the process on $\mathbb{Z}$ survives if and only if $\lambda>1/2$. Note that for $\lambda\leq 1/2$ the process restricted to $I_N$ will die out for all $N\geq 1$ since (by our construction)
the restricted process has less individuals than the unrestricted one.

Let $A_t(N)$ be the total number of individuals alive at time $t$ for the brw restricted to $I_N$. Define the critical value by
$$\lambda_c(N)=\inf\{\lambda>0:P_\lambda(A_t(N)\geq 1\mbox{ for all t})>0\}.$$

Since $P_\lambda(A_t(N)\geq 1\mbox{ for all t})$ is an increasing function of $\lambda$, if $\lambda>\lambda_c(N)$ the population restricted to $I_N$ survives forever with positive probability while the population dies out if $\lambda<\lambda_c(N)$. Actually, the population also dies out when 
$\lambda=\lambda_c(N)$ as will be explained below.

Tom Liggett (\cite{Liggett}) considered branching random walks on finite homogeneous trees. In the infinite tree each site has $d+1$ neighbors. The finite tree $T_{d,N}$ is obtained by retaining all those sites which can be reached from the origin (i.e. a fixed site on the tree) with a path of length less than or equal to $N$. In the particular case $d=1$, $T_{d,N}$ is exactly $I_N$ and the root is $0$. 

The analysis in  \cite{Liggett} (see page 319) is based on the representation of the branching random walk as a (non spatial) multi-type  branching process. The type of the individual is defined as its distance from the root. Thus, possible types are $0, 1, \dots,N$. Let $M_N(t)$ be the $(N+1)\times (N+1)$ matrix whose $(j,k)$ entry gives the mean number of type $k$ individuals at time $t$ produced by a single type $j$ individual at time $0$. The generator of this matrix generator is denoted by $A_N$. That is,
$$M_N(t)=\exp(A_N t).$$
The entries of the matrix $A_N$ are given by
$$
a(j,k)=\begin{cases}
    -1 &\mbox{ if }j=k\\
    \lambda & \mbox{ if }j=k+1\\
    \lambda & \mbox{ if }j=k-1\geq 1\\
    2\lambda  & \mbox{ if }j=0\mbox { and }k= 1\\
\end{cases}
$$
Let $f_N(x)=\det(xI-A_N)$ be the characteristic polynomial of the matrix $A_N$. It is shown in \cite{Liggett} that
for $N\geq 2$,
$$f_N(x)=(x+1)f_{N-1}(x)-\lambda^2f_{N-2}(x),$$
with $f_0(x)=x+1$ and $f_1(x)=(x+1)^2-2\lambda^2$.

Recall that a multi-type branching process survives if and only if the matrix
$A_N$ has a strictly positive eigenvalue, see for instance \cite{Athreya} (page 203). Note that we can write $A_N=-I+\lambda C_N$ where the entries of the matrix $C_N$ are constants. This shows that the eigenvalues of $A_N$ can be written as $-1+\lambda \sigma$ where $\sigma$ is an eigenvalue of $C_N$. Hence, $A_N$ has a positive eigenvalue if $C_N$ has a positive eigenvalue $\sigma$ and if $\lambda>0$ is such that $-1+\lambda \sigma>0$. Therefore, the critical value $\lambda_c(N)$ is the smallest value of $\lambda$ such that $A_N$ has zero as an eigenvalue.

 For instance, we see that $f_1(0)=0$ if $\lambda=\sqrt 2/2$. Hence,  $\lambda_c(1)=\sqrt 2/2$. Similarly, $\lambda_c(2)=\sqrt 3/3$. Moreover, since $\lambda_c(N)$ is the smallest value of $\lambda$ so that $A_N$ has zero as an eigenvalue we see that the process dies out for $\lambda=\lambda_c(N)$ for all $N\geq 1$.

In \cite{Liggett} it is also proved that the critical value of the branching random walk on $T_{d,N}$ has the following limit,
$$\lim_{N\to\infty}N^2\left (2\sqrt d\lambda_c(N)-1\right)=\frac{\pi^2}{2}.$$
We are interested in the case $d=1$. Not surprisingly, $\lambda_c(N)$ converges to $1/2$ (i.e. the critical value of the branching random walk on the whole line $\mathbb{Z}$). We also note that $\lambda_c(N)$ is strictly larger than 1/2 for all $N$. See also \cite{Mountford} for an analogous (but less precise) result on ${\mathbb Z}^d$ for all $d\geq 1$.

Let $1/2<\lambda\leq \sqrt 2/2$. 
Define
$$N_c=\min\{N\in\mathbb{N}:\lambda>\lambda_c(N)\}.$$
Since $\lambda_c(N)$ converges to $1/2$ the set $\{N\in\mathbb{N}:\lambda>\lambda_c(N)\}$ is not empty. By the well-ordering Principle it has a minimum $N_c$. Thus, the brw restricted to $I_N$ with birth rate $\lambda$ survives for all $N\geq N_c$.

Since $\lambda\leq \sqrt 2/2=\lambda_c(1)$ we see that $N_c\geq 2$. Therefore, $N_c-1$ is a natural number such that $\lambda\leq \lambda_c(N_c-1)$. Hence, the brw restricted to $I_N$ with birth rate $\lambda$ dies out for all $N\leq N_c-1$.
 This completes the proof of our result.

\bibliographystyle{amsplain}

\begin{thebibliography}{100}

\bibitem{Athreya} K.B. Athreya and P.E.Ney (1972) Branching processes. Springer-Verlag, New-York.

\bibitem{Haddad} N. M. Haddad, L. A. Brudvig, J. Clobert, K. F. Davies, A. Gonzalez
R. D. Holt, T. E. Lovejoy, J. O. Sexton, M. P. Austin, C. D. Collins, W. M. Cook, E. I. Damschen, R. M. Ewers, B. L. Foster, C. N. Jenkins, A. J. King, W. F. Laurance, D. J. Levey, C. R. Margules,
B. A. Melbourne, A. O. Nicholls, J. L. Orrock, D.X. Song, J. R. Townshend (2015)
Habitat fragmentation and its lasting impact on Earth’s ecosystems. Science Advances 1:e1500052.

\bibitem{Mountford} T. Mountford and R.B.Schinazi (2005) A note on branching random walks on finite sets. Journal of Applied Probability, 42, 287-294.

\bibitem{Liggett} T.M. Liggett (1999) Branching random walks on finite trees. In Perplexing problems in probability: papers in honor of Harry Kesten (M.Bramson and R.Durrett, eds), 315-330. Birkhauser, Boston.

\bibitem{Schinazi} R.B. Schinazi (2024) Classical and spatial stochastic processes (third edition). Birkhauser.

\end{thebibliography}

\end{document}